%% file: famdiv.tex
\documentclass[a4paper,12pt]{amsart}

\usepackage{xcolor}
\usepackage{amsmath}
\usepackage{amsfonts}
\usepackage{amssymb}
\usepackage{amsthm}
\usepackage{amsopn}
\usepackage{amscd}
\usepackage{pstricks}
\usepackage[utf8]{inputenc}
\usepackage{courier}
\usepackage{url}
\usepackage{color}
\usepackage{subfigure}
\usepackage{xypic}
\usepackage{graphicx}
\usepackage{paralist}
\usepackage[all,cmtip]{xy}
\usepackage{tikz}
\usepackage{hyperref}

\usepackage[vcentering,dvips]{geometry}
\entrymodifiers={+!!<0pt,\fontdimen22\textfont2>}
\title{Families of Invariant Divisors on Rational Complexity-One $T$-Varieties}

\author[A. Hochenegger]{Andreas Hochenegger}
\address{Mathematisches Insitut,
         Universit\"at zu K\"oln, 
         Weyertal 86--90,
         50931 K\"oln, Germany}
\email{ahochene@math.uni-koeln.de}

\author[N.O.~Ilten]{Nathan Owen Ilten}
\address{Department of Mathematics,
        University of California,
	Berkeley, CA 94720,
        USA}
\email{nilten@math.berkeley.edu}
\newcommand{\CC}{\mathbb{C}}
\newcommand{\QQ}{\mathbb{Q}}

\newcommand{\ZZ}{\mathbb Z}
\newcommand{\NN}{\mathbb N}
\renewcommand{\AA}{\mathbb{A}}
\newcommand{\PP}{\mathbb{P}}

\newcommand{\base}{\mathcal B}

\newcommand{\D}{\mathcal D}

\newcommand{\C}{\mathcal C}

\newcommand{\Y}{\mathcal Y}
\newcommand{\cX}{\mathcal{X}}
\newcommand{\mcP}{\mathcal P}
\newcommand{\CO}{\mathcal O}
\newcommand{\mcL}{\mathcal L}
\newcommand{\F}{\mathcal F}

\newcommand{\dfan}{\mathcal{S}}

\newcommand{\xvert}{\mathbf{vert}}
\newcommand{\xray}{\mathbf{ray}}
\newcommand{\lin}{\mathbf{lin}}

\newcommand{\tot}{\mathbf{tot}}
\newcommand{\pim}{\pi^\circ}
\newcommand{\bpim}{\bar\pi^\circ}

\newcommand{\vcent}[1]{\begin{array}{c}#1\end{array}}

\DeclareMathOperator{\Pic}{Pic}

\DeclareMathOperator{\Spec}{Spec}
\DeclareMathOperator{\Proj}{Proj}

\DeclareMathOperator{\WDiv}{WDiv}
\DeclareMathOperator{\CDiv}{CaDiv}
\DeclareMathOperator{\Div}{div}

\DeclareMathOperator{\spec}{Spec}
\DeclareMathOperator{\ord}{ord}
\DeclareMathOperator{\conv}{conv}
\DeclareMathOperator{\tail}{tail}

\DeclareMathOperator{\pos}{pos}

\DeclareMathOperator{\loc}{loc}

\DeclareMathOperator{\pr}{pr}

\DeclareMathOperator{\CaSF}{CaSF}
\DeclareMathOperator{\dCDiv}{-CaDiv}
\newcommand{\TCDiv}{T\hspace{-.15cm}\dCDiv}
\renewcommand{\sf}{\mathbf{sf}}

\theoremstyle{theorem}
\newtheorem{theorem}{Theorem}[section]
\newtheorem*{theorem*}{Theorem}

\newtheorem{prop}[theorem]{Proposition}
\newtheorem{lemma}[theorem]{Lemma}
\theoremstyle{definition}
\newtheorem{definition}[theorem]{Definition}
\newtheorem{rem}[theorem]{Remark}
\newtheorem{ex}[theorem]{Example}
\newgray{gray1}{.7}
\newgray{gray2}{.74}
\newgray{gray3}{.78}
\newgray{gray4}{.82}
\newgray{gray5}{.86}
\newgray{gray6}{.9}
\newgray{gray7}{.94}
\newgray{gray8}{.98}

\newgray{newgray1}{.4}
\newgray{newgray2}{.45}
\newgray{newgray3}{.55}
\newgray{newgray4}{.6}
\newgray{newgray5}{.65}
\newgray{newgray6}{.7}
\newgray{newgray7}{.75}
\newgray{newgray8}{.8}
\newgray{newgray9}{.85}
\newgray{newgray10}{.9}
\newgray{newgray11}{.95}
\newgray{newgray12}{1}

\include{figures}

\begin{document}
\maketitle
\begin{abstract}
We study invariant divisors on the total spaces of the homogeneous deformations of rational complexity-one $T$-varieties constructed by Ilten and Vollmert \cite{ilten:09b}.  In particular, we identify a natural subgroup of the Picard group for any  general fiber of such a deformation, which  is canonically isomorphic to the Picard group of the special fiber. This isomorphism preserves Euler characteristic, intersection numbers, and the canonical class.
\end{abstract}

\section*{Introduction}
Let $X$ be a rational complexity-one $T$-variety, that is, a normal rational variety $X$ endowed with an effective action by some torus $T$ of dimension one less than $\dim X$. Such varieties are the simplest examples of $T$-varieties after toric varieties, and can be described completely in terms of convex polyhedral data attached to points on the projective line, see \cite{tsurvey}. Any toric variety may be considered as belonging to the above class of varieties by restricting the action of the big torus to that of a codimension-one subtorus.

In \cite{ilten:09b}, R. Vollmert and the second author have shown how certain decompositions of the combinatorial data describing such an $X$ can be used to construct a deformation $\pi:\cX\to\base$ of $X$ with the action of $T$ on $X$ extending to an action on $\cX$ which preserves the fibers of $\pi$. We will call such deformations $T$-deformations. They have shown that if $X$ is a smooth complete toric variety, then the $T$-deformations of $X$ in fact span all directions in the vector space of first-order deformations of $X$. A similar result has also been attained by A. Mavlyutov, see \cite{mavlyutov:11a}.

On the other hand, invariant Cartier divisors on complexity-one $T$-varieties have been described in combinatorial terms by L. Petersen and H. S\"u\ss{}, see \cite{petersen:08a}. 
Fix now some one-parameter $T$-deformation $\pi:\cX\to \base$ of $X$. The purpose of this present paper is to explicitly compare the Picard group of a general fiber $\cX_s$ of $\pi$ with the Picard group of the special fiber $\cX_0=X$. To do this, we explicitly define a subgroup
$\TCDiv'(\cX_s)$ of the $T$-invariant Cartier divisors of $\cX_s$, see Definition \ref{def:casf}. This subgroup has the special property that any one of its elements can be canonically lifted to an invariant Cartier divisor on the total space $\cX$. By restricting the divisor on the total space to the special fiber $\cX_0$, we get a map $\pim$ from our subgroup $\TCDiv'(\cX_s)$ to the group $\TCDiv(\cX_0)$ of invariant divisors on $\cX_0=X$. This map can be described completely in combinatorial terms and respects linear equivalence. It thus descends to a map $\bpim:\Pic'(\cX_s)\to\Pic(\cX_0)$, where $\Pic'(\cX_s)$ is the image of $\TCDiv'(\cX_s)$ in $\Pic(\cX_s)$. Our main results can then be summed up in the following theorem:
\begin{theorem*}
Consider a one-parameter $T$-deformation $\pi:\cX\to\base$ of some projective rational complexity-one $T$-variety $X$. Then the corresponding map $$\bpim:\Pic'(\cX_s)\to\Pic(\cX_0)$$
is an isomorphism which preserves Euler characteristic, intersection numbers, canonical classes, and semiampleness.
\end{theorem*}       

Our original motivation for studying the behavior of divisors under $T$-deformation was to better understand the relationship between toric degenerations and exceptional sequences of line bundles on rational surfaces. In \cite{ilten:11e}, we apply our present results to give a sufficient criterion for a full exceptional sequence of line bundles on a rational $\mathbb{C}^*$-surface to remain exceptional under an equivariant degeneration. This may be used to show that for toric surfaces of Picard rank at most $4$, any exceptional sequence of line bundles of maximal length must be full.

The present paper is organized as follows.  In Section \ref{sec:prelim}, we introduce necessary notions from polyhedral geometry and prove a statement concerning the Minkowski summands of a polyhedral complex.  In Section \ref{sec:tvar}, we recall the basics  of $T$-varieties, including invariant divisors on complexity-one $T$-varieties. We then present a summary of the construction of $T$-deformations in a special case in Section \ref{sec:tdef}, and give a criterion for their properness. In Section \ref{sec:div1}, we introduce the subgroup  $\TCDiv'(\cX_s)$ and show how to lift its elements to the total space. We then proceed to analyze the map $\bpim$ in Section \ref{sec:div2}, completing the proof of the above theorem.  Finally, in Section \ref{sec:rems}, we make some closing remarks and present additional examples.

\section{Convex Geometry, Subdivisions, and Minkowski Decompositions}\label{sec:prelim}
In this section, we introduce notation and recall notions from convex geometry, as well as proving a result concerning Minkowski decompositions of polyhedral subdivisions. We refer to  \cite{ziegler} for basics on polyhedra.
Let $V$ be some finite-dimensional vector space over $\QQ$. For any polyhedron $\Delta\subset V$, its \emph{tailcone} is the set $\tail(\Delta):=\{v\in V\ |\ v+\Delta\subset\Delta\}$.
We denote its positive hull by $\pos(\Delta)$. For a cone $\sigma\subset V$, we denote its dual by $\sigma^\vee$.
Consider some finite set $S$ of polyhedra in $V$.
The \emph{support} of $S$ is the set
$$
|S|:=\bigcup_{\Delta\in S} \Delta.
$$
The set $S$ is a \emph{polyhedral subdivision} if it is closed under taking faces, and any two elements of $S$ intersect in a common face. 
We say that some set $S$ of polyhedra \emph{induces} a subdivision if the union of all faces of elements of $S$ is a subdivision; we denote the resulting subdivision by $\langle S \rangle$. 
We will be particularly interested in subdivisions with convex support.

Consider now any subdivision $S$ in $V$ with $|S|$ convex. A \emph{support function} on $S$ is a continuous map $f:|S|\to \QQ$ such that the restriction of $f$ to any $\Delta\in S$ is affine.
The \emph{linear part} of $f$ is the map $f^\lin:\tail(|S|)\to \QQ$ defined by 
$$
f^\lin(v):=\lim_{\lambda\to\infty} f(v_0+\lambda v)/\lambda
$$
for any $v_0\in |S|$.
Consider a lattice $N$ with dual $M$ and suppose that $V=N_\QQ=N\otimes\QQ$. We say that $f$ is \emph{integral} (with respect to $N$) if it has integral slope and translation, i.e. for any $\Delta\in S$, $f$ restricted to $\Delta$ can be written as $\langle \cdot , u \rangle+a$ for some $u\in M$ and $a\in \ZZ$. 

Consider now some finite set $S$ of full-dimensional polyhedra in $V$. In general, 
the set $\conv |S|$ may no longer be closed (and thus not a polyhedron). By a \emph{facet} of $\conv |S|$, we mean any facet of $\overline{\conv |S|}$ which is contained in $\conv |S|$.
We will use the following result to help recognize when a collection of polyhedra induces a subdivision.
\begin{lemma}[cf. {\cite[Lemma 2.1]{huber:00a}}]\label{lemma:issub}
Let $S$ be some finite set of full-dimensional polyhedra in $V$. Then 
$\langle S \rangle$ is a subdivision with convex support if and only if the following conditions are satisfied:
\begin{enumerate}
\item There is a point in $|S|$ which is contained in exactly one $\Delta\in S$.
\item For every $\Delta\in S$ and every facet $\tau\prec \Delta$, either $\tau$ is contained in a facet of $\conv |S|$ or there is another $\nabla\in S$ with $\tau\prec\nabla$. We say in this case that $\Delta$ and $\nabla$ are \emph{adjacent}.
\item If $\Delta,\nabla\in S$ are adjacent then they lie in opposite halfspaces with respect to their common facet.
\end{enumerate}

\end{lemma}

\begin{proof}
It is immediate that any $S$ inducing a subdivision with convex support satisfies the desired conditions. 
Conversely, let $S$ be any finite set $S$ of full-dimensional polyhedra in $V$ satisfying the assumptions of the lemma.
If the elements of $S$ are bounded, then the claim that $\langle S \rangle$ is a subdivision with convex support is just \cite[Lemma 2.1]{huber:00a}. In order to show that claim for general $S$, we will reduce to this special case.

For any full-dimensional polytope $\Box\subset V$ containing the origin, define
$$S\cap\Box=\{\Delta\cap \Box\ |\ \Delta\in S\}.$$
For $k\in\NN$  sufficiently large, $S$ will induce a subdivision if and only if $S\cap k\cdot \Box$ induces a subdivision. We thus must show that $S\cap k\cdot \Box$ satisfies the conditions of the lemma for sufficiently large $k$.

Conditions (1) and (3) are immediate for sufficiently large $k$. Now, for sufficiently large $k$ and any $\Delta\in S$, $\Delta \cap k\cdot \Box$ has two kinds of facets: those of the form $\tau\cap k\cdot\Box$ for $\tau$ a facet of $\Delta$, and those of the form $\Delta\cap k\cdot\tau'$ for $\tau'$ a facet of $\Box$. A facet of the form $\tau\cap k\cdot\Box$ satisfies condition (2), since $\tau$ satisfies (2) for the set $S$. On the other hand, a facet of the form $\Delta\cap k\cdot\tau'$ is clearly contained in a facet of $\conv |S\cap k\cdot\Box|$. Thus, condition (2) is fulfilled for the set $S\cap k\cdot \Box$.
\end{proof}

Recall that an $r$-term \emph{Minkowski decomposition} of a polyhedron $\Delta\subset V$ consists of polyhedra $\Delta^0,\ldots,\Delta^r$ all with tailcone $\tail(\Delta)$ such that $\Delta=\Delta^0+\ldots+\Delta^r$. Suppose that $V=N_\QQ$ for some lattice $N$ with dual $M$. Then $\Delta$ is \emph{admissible} (with respect to $N$) if for all $u\in M\cap\tail(\Delta)^\vee$, the  minimum value of $u$ on $\Delta^i$ is integral for all $0\leq i \leq r$ with at most one exception.
Now consider any polyhedral subdivision $S$. 
\begin{definition}
An $r$-term \emph{Minkowski decomposition}\footnote{This is weaker than the definition of a Minkowski decomposition in \cite{ilten:09b}. In \cite{ilten:10b}, this is called a \emph{predecomposition}.} of $S$ consists of Minkowksi decompositions $\Delta=\Delta^0+\ldots+\Delta^r$ for every $\Delta\in S$ such that for any $\Delta,\nabla\in S$ with $\Delta\cap\nabla \neq \emptyset$, we have
$$
\Delta^i\cap \nabla^i=(\Delta\cap\nabla)^i
$$
for all $0\leq i \leq r$. Such a decomposition is \emph{admissible} if the corresponding decomposition of every $\Delta\in S$ is admissible.
An admissible decomposition is \emph{essentially locally trivial} if for all $\Delta\in S$, some $\Delta^i$ is a lattice translate of $\Delta$.
\end{definition}

\begin{figure}[htbp]
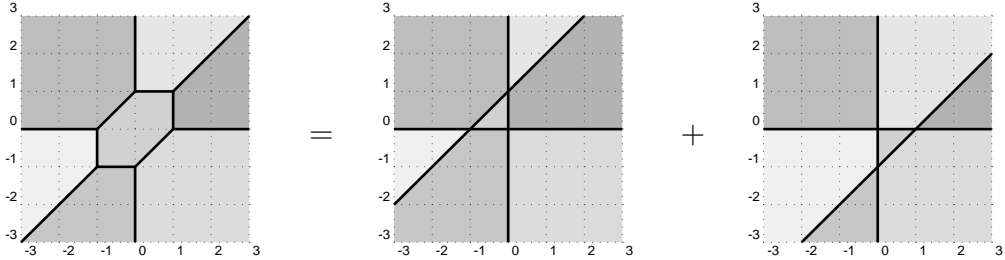

    \centering
\begin{displaymath}
    \vcent{{\dpnull}}=\vcent{{\dpnullia}}+\vcent{{\dpnullib}}
\end{displaymath}  
  \caption{A Minkowski decomposition}\label{fig:md}
\end{figure}

\begin{ex}[A Minkowski decomposition]\label{ex:md}
In Figure \ref{fig:md} we picture a Minkowski decomposition of a polyhedral subdivision $S$ in $\QQ^2$. For any full-dimensional polyhedron $\Delta$, the corresponding summands $\Delta^0$ and $\Delta^1$ have the same shade of gray. For example, the hexagon in the middle of the subdivision decomposes into the sum of two triangles. Since all polyhedra involved are lattice polyhedra, this decomposition is automatically admissible.
\end{ex}

The following proposition uses Lemma \ref{lemma:issub} coupled with the Cayley trick, cf. \cite[Proposition 3.5]{huber:00a} and \cite[Lemma 4.6]{mavlyutov:11a}. We include a proof for the sake of completeness.
\begin{prop}\label{prop:mdcomp}
Let $S$ be a polyhedral subdivision with  convex support.
Consider any $r$-term Minkowski decomposition of $S$. Then for any $0\leq i \leq r$, the set
$$S^i=\{\Delta^i\ |\ \Delta\in S\}$$
is a polyhedral subdivision with convex support.
\end{prop}
\begin{proof}
We begin with some simplifying assumptions.
We may assume that $r=1$ and $i=0$. Indeed, the data $\Delta=\Delta^i+(\sum_{j\neq i} \Delta^j)$ for all $\Delta\in S$ provides a Minkowski decomposition of $S$. Furthermore, we may assume that $|S|$ is full-dimensional, since if $|S|$ is contained in some subspace $V'$ of $V$, then the $S^i$ are contained in affine translates of $V'$. We may also assume that $|S|$ contains the origin in its interior. Finally, we may assume that $|S|=V$. If this is not the case, we can consider the subdivision $S'$ induced by the polyhedra of $S$ together with the polyhedra of the form $\Delta+\pos(\Delta)$ for any $\Delta\in S$ contained in the boundary of $|S|$. Any Minkowski decomposition of $S$ induces a Minkowski decomposition of $S'$, and $S^i$ is a subset of $(S')^i$.

For any polyhedra $\Delta^0,\Delta^1\subset V$ with common tailcone,\footnote{Note that if $\tail(\Delta^0)\neq\tail(\Delta^1)$, then the set $\C(\Delta^0,\Delta^1)$ will not be closed.}  their \emph{Cayley polyhedron} is
$$
\C(\Delta^0,\Delta^1):=\conv \{(\Delta_0,e_0)\cup(\Delta_1,e_1)\}\subset V\oplus \QQ^2
$$
where $e_0,e_1$ is the standard basis of $\QQ^2$ with dual basis $e_0^*,e_1^*$. The affine embedding $V \to V\oplus \QQ^2$ which sends any $v$ to $(v+e_0+e_1)/2$ identifies $\Delta^0+\Delta^1$ with $\C(\Delta^0,\Delta^1)\cap W$, where $W=\{e_0^*=e_1^*=1/2\}$.
Now let $T$ be the subset of $S$ consisting of all full-dimensional polyhedra. We will show that the set
$$
\C(T):=\{C(\Delta^0,\Delta^1)\ |\ \Delta\in T\}
$$
induces a polyhedral subdivision in $V\oplus \QQ^2$. It then follows that $S^0$ is a polyhedral subdivision. Indeed, $S^0$ can be identified with $\langle\C(T)\rangle \cap \{e_0^*=1\}$.

In the following, we will identify $S$ with  $\langle \C(T)\cap W\rangle$.
It is clear that any $\Box\in \C(T)$ is full-dimensional. Thus we can apply the criteria from Lemma \ref{lemma:issub}. Property (1) is immediate: Take any $x\in |S|$ contained in exactly one $\Delta\in T$. Then $x$ is only contained in $\C(\Delta^0,\Delta^1)\in\C(T)$.  To see property (3), suppose $\C(\Delta^0,\Delta^1)$ and $\C(\nabla^0,\nabla^1)$ are adjacent along $\tau$ for some $\Delta,\nabla\in T$. Then $\Delta$ and $\nabla$ are adjacent along $\tau\cap W$, and lie on opposite sides of the corresponding hyperplane if and only if $\C(\Delta^0,\Delta^1)$ and $\C(\nabla^0,\nabla^1)$ do.

Finally, to check (2), let $\tau$ be any facet of some $\C(\Delta^0,\Delta^1)\in\C(T)$.
Then either $\tau$ is contained in $\{e_0^*=1\}$ or $\{e_1^*=1\}$, or $\tau\cap W$ is a facet in $S$. In the former case, $\tau$ is contained in a facet of $\conv |\C(T)|$.  Otherwise, $\Delta$ is adjacent to some $\nabla$ along $\tau\cap W$, and it follows from the definition of Minkowski decomposition that $\C(\Delta^0,\Delta^1)$ is adjacent to $\C(\nabla^0,\nabla^1)$ along $\tau$.
\end{proof}

\begin{rem}
Let $S$ be any polyhedral subdivision with convex support, and consider the Minkowski decomposition given by $\Delta=\Delta^0+\Delta^1=\Delta+\tail(\Delta)$. Then by the above proposition,
$$\tail(S):=\{\tail(\Delta)\ |\ \Delta\in S\}$$
is a polyhedral subdivision of $\tail(|S|)$, the so-called \emph{tailfan} of $S$. The fact that this set of cones forms a subdivision is a special case of \cite[Theorem 3.4]{burgos:10a}.
\end{rem}

\section{$T$-Varieties}\label{sec:tvar}
We now recall the construction of $T$-varieties from p-divisors and divisorial fans, see \cite{altmann:06a} and \cite{altmann:08a}, as well as recalling the description of invariant Cartier divisors on complexity-one $T$-varieties \cite{petersen:08a}. For an introduction to and a  survey of the theory of $T$-varieties, see \cite{tsurvey}.

Fix some lattice $N$ with dual $M$.  We will be considering the algebraic torus $T=\CC^*\otimes N=\spec \CC[M]$. Recall that a $T$-variety of complexity $k$ is a normal variety $X$ together with an effective $T$-action such that $\dim X-\dim T=k$.
Now let $Y$ be a normal projective variety over $\CC$, and $\sigma\subseteq N_\QQ$ a pointed polyhedral cone.
\begin{definition}
A \emph{polyhedral divisor} on $Y$ with tailcone $\sigma$ is a formal sum
$$
\D=\sum \D_P\otimes P
$$
over all prime divisors $P\subset Y$, such that any $\D_P$ is either the empty set or a polyhedron in $N_\QQ$ with tailcone $\sigma$, and $D_P\neq \sigma$ for only finitely many $P$. The \emph{locus} of $\D$ is the set 
$$
\loc(\D):=Y\setminus \bigcup_{\substack {P\subset Y\\ \D_P=\emptyset}} P.
$$
\end{definition}

Any polyhedral divisor $\D$ gives rise to an evaluation map $\D:\sigma^\vee\to \WDiv_\QQ(\loc(\D))$ by setting
$$
\D(u):=\sum_{P\subset \loc(\D)} \min_{v\in \D_P}\langle v,u\rangle P.
$$
Here, $\WDiv_\QQ(\loc(\D))$ denotes the group of Weil divisors on $\loc(\D)$ with rational coefficients.
\begin{definition}
A polyhedral divisor $\D$ is a \emph{p-divisor} if $\D(u)$ is $\QQ$-Cartier and semiample for all $u\in\sigma^\vee$, $\D(u)$ is big for some $u\in \sigma$, and $\loc(\D)$ is semiprojective.\footnote{Recall that a divisor is big if some multiple admits a section with affine complement. Likewise, a divisor is semiample if some multiple is base-point free. A variety is semiprojective if it is projective over an affine variety.}
\end{definition}

\begin{rem}
If $Y$ is a curve, we define the degree of a polyhedral divisor $\D$ to be $$\deg \D=\sum_{P\in Y} \D_P.$$ Then $\D$ is a p-divisor if and only if $\deg \D \subsetneq \sigma$, and for any $u\in\sigma^\vee$ with $\deg \D(u)=0$, $\D(u)$ has a principal multiple.
\end{rem}

Any p-divisor $\D$ gives rise to an affine $T$-variety of complexity equal to the dimension of $Y$, see \cite[Theorem 3.1]{altmann:06a}:
$$
X(\D):=\spec \bigoplus_{u\in\sigma^\vee\cap M} H^0(\loc(\D),\CO(\D(u))). 
$$
Conversely, any affine $T$-variety can be constructed from some p-divisor \cite[Theorem 3.5]{altmann:08a}. 

To describe non-affine $T$-varieties, we need a little more notation. Consider two polyhedral divisors $\D,\D'$ on $Y$. Their intersection $\D\cap\D'$ is gotten by intersecting their coefficients, that is,
$$
\D\cap\D':=\sum \D_P\cap\D_P'\otimes P.
$$
Likewise, we say $\D'\subseteq\D$ if $\D_P'\subseteq\D_P$ for all prime divisors $P$. If $\D,\D'$ are in fact p-divisors and $\D'\subseteq\D$, then there is a dominant morphism
$$
X(\D')\to X(\D)
$$
induced by the opposite inclusion of coordinate rings.  We say that $\D'$ is a \emph{face} of $\D$, written $\D'\prec \D$, if this map is an open embedding. A necessary condition for $\D'\prec\D$ is that $\D_P'\prec\D_P$ for all prime divisors $P$.
If $Y$ is a curve, then $\D'\prec\D$ if and only if the previous condition is satisfied, and $\deg \D'=\deg \D\cap \tail (\D')$.
For a complete characterization in the general case see \cite[Definition 5.1]{altmann:08a}.

\begin{definition}
A \emph{divisorial fan} $\dfan$ on $Y$ is a finite set of p-divisors on $Y$ closed under taking intersections, such that any two elements of $\dfan$ intersect in a common face. For any prime divisor $P\subset Y$, the slice of $\dfan$ at $P$ is the set $\dfan_P:=\{\D_P\ |\ \D\in\dfan\}$.
The tailfan of $\dfan$ is the set $\tail(\dfan):=\{\tail(\D)\ |\ \D\in\dfan\}$. 
\end{definition} 
\begin{rem}
The face condition ensures that the slices and tailfan of a divisorial fan are polyhedral subdivisions.
\end{rem}
To any divisorial fan $\dfan$, we can associate a normal $T$-scheme
$$
X(\dfan)=	\bigcup_{\D\in\dfan} X(\D)/\sim,
$$
where $\sim$ is the gluing along open embeddings of the sort
$X(\D)\hookleftarrow X(\D'\cap\D)\hookrightarrow X(\D)$. Separatedness of the scheme $X(\dfan)$ can be characterized in terms of $\dfan$, see \cite[Section 7]{altmann:08a}.
On the other hand, every $T$-variety can be described via some divisorial fan \cite[Theorem 5.6]{altmann:08a}.

\begin{ex}[A compactified cone over the del Pezzo surface of degree six]\label{ex:dp}
	Let $\dfan$ consist of p-divisors of the form $\Delta\otimes\{0\}$ on $\PP^1$ for full-dimensional unbounded $\Delta$ in the polyhedral subdivision $S$ of Example \ref{ex:md} and Figure \ref{fig:md}, the p-divisor $\nabla\otimes\{0\}+\emptyset\otimes\infty$ for $\nabla$ the compact hexagon in $S$, and intersections thereof.  Then $\dfan$ is a divisorial fan, $S=\dfan_0$, and  $\overline{C(dP_6)}:=X(\dfan)$ is a compactification of the anticanonical cone over the del Pezzo surface of degree six.
\end{ex}

We now specialize to the case of complexity-one $T$-varieties. Thus, we will be considering divisorial fans on smooth projective curves. Let $\dfan$ be such a divisorial fan. Then $X(\dfan)$ is always separated, and it is complete if and only if $|\dfan_P|=N_\QQ$ for all prime divisors $P$, see \cite[Remark 7.4]{altmann:08a}. We call such fans complete.  Suppose now that $|\dfan_P|$ is convex for all $P$. We now recall the description of invariant Cartier divisors on $X(\dfan)$ from \cite{petersen:08a}.

\begin{definition}
A \emph{divisorial support function} on $\dfan$  is a formal sum
$$h=\sum_{P\subset Y} h_P\otimes P$$
such that $h_P$ is a support function of $\dfan_P$, $h_P^\lin$ is independent of $P$, and $h_P$ differs from $h_P^\lin$ for only finitely many points $P$. A divisorial support function $h$ is \emph{integral} if $h_P$ is integral for all points $P$ in $Y$. 
An integral divisorial support function $h$ is \emph{Cartier} if for any $\D\in\dfan$ with complete locus, there exists $f\in\CC(Y)^*$ and $u\in M$ such that $(h_P)_{|\D_P}\equiv \langle \cdot, u \rangle + \ord_P f$ for all $P\in Y$.
\end{definition}

Let $\CaSF(\dfan)$ denote the group of Cartier integral divisorial support functions on $\dfan$ (with pointwise addition). To any pair $(u,f)\in M\times \CC(Y)^*$ we can associate such a support function $\sf(u,f)\in\CaSF(\dfan)$, where the $P$-coefficient of $\sf(u,f)$ is given by $\langle \cdot,u\rangle +\ord_P f$. We call such support functions \emph{principal}.
The group $\CaSF(\dfan)$ is in fact isomorphic to $\TCDiv(X(\dfan))$, the group of $T$-invariant Cartier divisors on $X(\dfan)$, with the isomorphism inducing a bijection between principal support functions and invariant principal divisors \cite[Proposition 3.10]{petersen:08a}. Given $h\in\CaSF(\dfan)$, we denote the corresponding Cartier divisor by $D_h$. 
For a principal invariant divisor $h=\sf(u,f)$, $D_h$ is simply $\Div(f\chi^u)$.
To construct $D_h$ for some general $h$, one must refine the divisorial fan $\dfan$ such that for all $\D\in\dfan$, the restriction of $h$ to $\D$ is principal. 

Sometimes it is more practical to deal with invariant prime (Weil) divisors. These fall into two classes, see \cite[Section 3.2]{petersen:08a}. Firstly, there are the ``vertical'' divisors arising as the closures of codimension-one $T$-orbits. These are parametrized by pairs $(P,v)$, where $P\in Y$ and $v\in \xvert_P(\dfan)$, i.e. $v$ is a vertex of $\dfan_P$. We denote the prime divisor corresponding to such a pair by $D_{P,v}$. Secondly, there are the ``horizontal'' divisors arising as closures of families of codimension-two $T$-orbits. These are parametrized by rays $\rho\in\xray(\dfan)$, where $\xray(\dfan)$ is the set of all rays $\rho\in\tail(\dfan)$ such that $\rho\cap\deg \dfan=\emptyset$, with $\deg \dfan$ defined as
$$\deg\dfan=\bigcup_{\D\in\dfan}\sum_{P\in Y} \D_P.$$
We denote the divisor corresponding to such a $\rho$ by $D_\rho$.
The following proposition relates divisorial support functions to these prime divisors:

\begin{prop}[{\cite[Corollary 3.19]{petersen:08a}}]\label{prop:weil}
Consider any $h\in\CaSF(\dfan)$. Then
$$
D_h=-\sum_{\rho\in\xray(\dfan)} h^\lin(v_\rho)D_\rho-\sum_{\substack{P\in Y\\v\in\xvert_P(\dfan)}} \mu(v)h_P(v)D_{P,v},
$$
where $\mu(v)$ is the smallest natural number such that $\mu(v)v\in N$ and $v_\rho$ is the primitive lattice generator for $\rho$.
\end{prop}

\section{$T$-Deformations}\label{sec:tdef}
We briefly recall a special case of the construction of $T$-deformations found in \cite{ilten:09b}, and then prove a result concerning the separatedness and properness of these deformations.
Starting now, $Y$ will always be $\PP^1=\Proj \CC[y_0,y_1]$. 

Consider a p-divisor $\D$ on $Y$ and the finite set
$$\mcP(\D)=\{P\in Y\ |\ \D_P\neq \tail(\D)\}\setminus\{0\}.$$
Let $\base_\D$ be the complement of $\mcP(\D)\setminus\{\infty\}$ in $\AA^1=\spec \CC[t]$, and take $\Y_\D=Y\times \base_\D$.
Now consider any admissible one-term Minkowski decomposition of $\D_0$ where $0=V(y_1)$. If $\D_0=\emptyset$, then we will consider the `decomposition' $\D_0=\emptyset+\emptyset$. 
Such a Minkowski decomposition will give rise to a deformation of $X(\D)$ over the base space $\base_\D$. Consider the polyhedral divisors
\begin{align*}
\D^\tot&=\D_0^0\otimes (\{0\}\times \base)+\D_0^1\otimes V(ty_0-y_1)+\sum_{P\in\mcP(\D)} \D_P\otimes (P\times \base)\\
\D^{(s)}&=\D_0^0\otimes \{0\}+\D_0^1\otimes s+\sum_{P\in\mcP(\D)} \D_P\otimes P
\end{align*}
where $\D^\tot$ is a polyhedral divisor on $\Y$, and $\D^{(s)}$ is a polyhedral divisor on $Y$ for any $s\in \base$. In particular, $\D^{(0)}=\D$. 
All these polyhedral divisors are actually p-divisors.
The $T$-variety $X(\D^\tot)$ comes with a regular map $\pi_\D:X(\D^\tot)\to \base_\D$ gotten via the composition of the map $X(\D^\tot)\dashrightarrow\Y_D$ with the projection $\Y_D\to \base_D$.
\begin{theorem}[{\cite[Theorem 2.8]{ilten:09b}}]
The family $\pi_\D:X(\D^\tot)\to\base_\D$ is flat, with $\pi_\D^{-1}(s)\cong X(\D^{(s)})$ for all $s\in\base_\D$. In particular, $\pi^{-1}(0)\cong X$.
\end{theorem}

These deformations can be glued together to give deformations of non-affine complexity-one rational $T$-varieties. Let $\dfan$ be a divisorial fan  on $Y$ such that $|\dfan_0|$ is convex.
Let $\mcP=\bigcup_{\D\in\dfan}\mcP(\D)$, and take $\base=\AA^1\setminus\mcP$. 
 Any admissible one-term Minkowski decomposition of $\dfan_0$ gives rise to a deformation of $X=X(\dfan)$ over $\base$. Indeed, the decomposition of $\dfan_0$ uniquely determines a decomposition of $\D_0$ for any $\D\in\dfan$. Now, let $\dfan^\tot$ be the set of p-divisors on $\Y=Y\times \base$ induced via intersection of the p-divisors $\D^\tot$ for any $\D\in\dfan$. Likewise, for any $s\in\base$, let
$\dfan^{(s)}$ be the set of p-divisors on $Y$ induced via intersection of the p-divisors $\D^{(s)}$ for any $\D\in\dfan$.
\begin{theorem}[{cf. \cite[Theorem 4.4]{ilten:09b}}]\label{thm:tdefprop}
 The sets $\dfan^\tot$ and $\dfan^{(s)}$ are divisorial fans. The deformations $\pi_\D$ glue together to give a flat family $\pi:X(\dfan^\tot)\to\base$ with $\pi^{-1}(s)\cong X(\dfan^{(s)})$ for all $s\in\base$. In particular, $\pi^{-1}(0)\cong X(\dfan)$.
   Furthermore, $X(\dfan^\tot)$ is separated, and $\pi$ is proper if and only if $\dfan$ is complete.     \end{theorem}

Note that the torus $T$ acts on the total space $X(\dfan^\tot)$ and the map $\pi$ is in fact $T$-invariant.
We thus call deformations of the above sort \emph{$T$-deformations}. We now show that in the situation presented here, $T$-deformations are always separated, and provide a criterion for properness:

\begin{proof}
  The claims regarding $\dfan^\tot$ and $\dfan^{(s)}$ and the fibers of $\pi$ are immediate from Section 4 of \cite{ilten:09b}, except that the definition of Minkowski decomposition used there appears stronger than our present definition.  However, in the special case of $|\dfan_0|$ convex, these are in fact equivalent due to Proposition \ref{prop:mdcomp}.

We now discuss separatedness and properness.
Let $\nu:\CC(\Y)^*\to\QQ$ be a valuation with center $y\in \Y$. 
This induces a group homomorphism 
$
\nu:\CDiv_{\geq 0} \Y\to\QQ_{\geq 0}
$
sending a divisor $D$ with local equation $f$ at $y$ to $\nu(f)$.
Then $\nu$ defines a set of polyhedra $\dfan^\tot_\nu$ called a weighted slice, see \cite[Section 7]{altmann:08a}:
\begin{align*}
\dfan_\nu^\tot:=\{\D_\nu\ |\ \D\in\dfan^\tot\},\qquad 
\D_\nu:=\sum_{P\subset \Y} \nu(P)\D_P.
\end{align*}
We claim that any such weighted slice $\dfan_\nu^\tot$ is a polyhedral subdivision; this will imply that $X(\dfan^\tot)$ is separated by the evaluation criterion of \cite[Section 7]{altmann:08a}. Indeed, if $y\neq 0=V(t,y_1)$, then $\dfan^\tot_\nu$ is simply a dilation of a slice $\dfan_P$ for some $P\in Y$. Suppose instead $y=0=V(t,y_1)$, and let $a_0=\nu(V(y_1))$ and $a_1=\nu(V(y_0t-y_1))$. Then $\dfan^\tot_\nu$ consists of  intersections of polyhedra of the form $a_0\Delta^0+a_1\Delta^1$ for $\Delta\in \dfan_0$. But for any $\Delta\in\dfan_0$, the decomposition $$(a_0+a_1)\Delta=(a_0\Delta^0+a_1\Delta^1)+(a_1\Delta^0+a_0\Delta^1)$$ gives a one-term Minkowski decomposition of the polyhedral subdivision   $(a_0+a_1)\dfan_0:=\{(a_0+a_1)\Delta\ |\ \Delta\in\dfan_0\}$.
Thus, by proposition \ref{prop:mdcomp}, $\dfan_\nu^\tot$ is a polyhedral subdivision.

        The claim regarding properness uses a relative version of the evaluation criterion of \cite[Section 7]{altmann:08a}; this is described in \cite[Theorem 46]{tsurvey}. Using the notation from \cite{tsurvey}, the map $\pi$ is in fact a torus equivariant morphism corresponding to the triple $(\pr,F,0)$, where $\pr:Y\times \base\to \base$ is the projection and $F:N\to 0$ is the zero map. Since $\pr$ is proper, $\pi$ is proper if and only if  each weighted slice of $\dfan^\tot$ is a complete polyhedral subdivision. By the above discussion, this is the case exactly when $\dfan$ is complete.
\end{proof}

\begin{rem}
By considering multi-term Minkowski decompositions of multiple slices of $\dfan$, one may construct multi-parameter deformations of $X(\dfan)$. Similar results hold for the separatedness and properness of these deformations.
We also remark that if a one-term Minkowski decomposition is essentially locally trivial, and $\dfan$ only has two non-trivial slices, then the resulting deformation is locally trivial after being pulled back to the fat point $\Spec \CC[t]/t^2$, see \cite[Theorem 5.1]{ilten:09b}. Note that although locally trivial, such a deformation is nonetheless in general not trivial.
\end{rem}

\begin{ex}[A deformation of $\overline{C(dP_6)}$]\label{ex:dpdef}
Let $\dfan$ be the divisorial fan from Example \ref{ex:dp}. The Minkowski decomposition of $\dfan_0$ from Example \ref{ex:md} pictured in Figure \ref{fig:md} gives a $T$-deformation $\pi$ of $\overline{C(dP_6)}=X(\dfan)$.
 For $s\neq 0$, the fiber $\pi^{(-1)}(s)$ is isomorphic to $\PP^1\times\PP^1\times\PP^1$. Indeed, the fiber is described by a divisorial fan $\dfan^{(s)}$ with exactly two non-trivial slices $\dfan_0^{(s)}$ and $\dfan_s^{(s)}$ which are simply the summands in the decomposition of Figure \ref{fig:md}.
In order to see that $X(\dfan^{(s)})=\PP^1\times\PP^1\times\PP^1$, one may reverse the downgrade procedure described in \cite[Section 5]{altmann:08a}.
	
The deformation $\pi$ was presented in \cite{jahnke:06a} as an example of a smoothing of a Fano variety with canonical singularities. A combinatorial description similar to the one presented here can be found in \cite{suess:08a}.
\end{ex}

\section{Invariant Families of Divisors I}\label{sec:div1}
Let $\dfan$ be a divisorial fan on $Y=\PP^1$ with $|\dfan_0|$ convex, and consider an admissible one-term Minkowski decomposition of $\dfan_0$ leading to a $T$-deformation $\pi:X(\dfan^\tot)\to\base$ as in the previous section. We denote the total space $X(\dfan^\tot)$ by $\cX$, and for any $s\in\base$, we write $\cX_s:=\pi^{-1}(s)$. Let $\base^*$ denote the complement of the origin in $\base$. 

Our goal is now to compare the Picard groups $\Pic (\cX_s)$ as $s\in \base$ varies. Our strategy is the following: for any fixed $s\in \base ^*$ we will identify a subgroup $\TCDiv'(\cX_s)$ of $\TCDiv(\cX_s)$ such that any element of $\TCDiv'(\cX_s)$ naturally lifts to an invariant Cartier divisor on the total space $\cX$. We can then restrict this divisor to the special fiber $\cX_0$, giving us an element of $\TCDiv(\cX_0)$. Thus, we will have a natural map $\pim:\TCDiv'(\cX_s)\to \TCDiv(\cX_0)$. Since this map respects linear equivalence, we can then use it to compare subgroups of the Picard groups of the fibers.

Our first task is now to identify the special subgroups $\TCDiv'$ which allow for natural lifting of divisors to $\cX$. This will be taken care of by the following definition:
\begin{definition}\label{def:casf}
        For $s\in \base^*\subset Y$, define $\CaSF'(\dfan^{(s)})$ to consist of those $h\in\CaSF(\dfan^{(s)})$ such that for all $\Delta \in \dfan_0$, we can find $u\in M$ and $a_0,a_1\in\ZZ$ satisfying
\begin{align*}
    (h_{P_i})_{|\Delta^i}(v)=\langle v,u \rangle +a_i\qquad i\in\{0,1\}
       \end{align*}
where $P_0=0$, $P_1=s$.
Note that $\CaSF'(\dfan^{(s)})=\CaSF(\dfan^{(s)})$ if the decomposition of $\dfan$ is essentially locally trivial. Finally, by $\TCDiv'$ we denote the image of $\CaSF'$ under the natural map described in Section \ref{sec:tvar}.
\end{definition}

\begin{rem} Let $\TCDiv'(\cX)$ denote the group of $T$-invariant Cartier divisors on $\cX$ which intersect the fibers of $\pi$ properly. Then it will follow from the following discussion that $\TCDiv'(\cX_s)$ is the image of $\TCDiv'(\cX)$ under restriction to $\cX_s$.  
\end{rem}

Fix now some $s\in \base^*$ and choose some support function $h\in \CaSF'(\dfan^{(s)})$; this corresponds to an invariant Cartier divisor $D_{h}\in\TCDiv'(\cX_s)$. We will be showing that this can be lifted to a Cartier divisor $D_h^\tot$ on $\cX$. We first will need invariant open coverings of $\cX_s$ and $\cX$. 
Take $\mcP$ as in the previous section.
For $P\in \mathcal{P}$ and $\D\in \dfan$ with noncomplete locus, set
\begin{align*}
        U_{\D,P}&=X\Big(\D^{(s)}+\emptyset\otimes \{0\}+\emptyset\otimes s +\sum_{\substack{Q\in \mathcal{P}\\Q\neq P}} \emptyset \otimes Q\Big)\\
U_{\D,P}^{\tot}&=X\Big(\D^{\tot}+\emptyset\otimes (\{0\}\times \base)+\emptyset\otimes V(ty_0-y_1)+\sum_{\substack{Q\in \mathcal{P}\\Q\neq P}} \emptyset \otimes (Q\times \base) \Big)
\end{align*}
and likewise set
\begin{align*}
U_{\D,0}&=X\Big(\D^{(s)}+\sum_{Q\in \mathcal{P}} \emptyset \otimes V(y_0-\lambda_Qy_1)\Big)\\
U_{\D,0}^{\tot}&=X\Big(\D^{\tot}+\sum_{Q\in \mathcal{P}} \emptyset \otimes V(y_0-\lambda_Qy_1)\Big).
\end{align*}
On the other hand, for $P\in\mathcal{P}\cup\{0\}$ and $\D\in \dfan$ with complete locus, set $U_{\D,P}=X(\D^{(s)})$ and $U_{\D,P}^{\tot}=X(\D^{\tot})$. One easily checks that $\{U_{\D,P}\}$ and $\{U_{\D,P}^{\tot}\}$ define invariant open coverings of respectively $\cX_s$ and $\cX$, and that $U_{\D,P}^\tot\cap\cX_s=U_{\D,P}$. These open coverings may in fact be finer than necessary for defining the desired Cartier divisor.

For each $P\in\mathcal{P}\cup\{0\}$ and $\D\in \dfan$, let $u_{\D,P}\in M$, $f_{\D,P}\in \CC(Y)$ be such that ${D_{h}}_{|U_{\D,P}}=\Div(f_{\D,P}\cdot \chi^{u_{\D,P}})$. Such $f_{\D,P},u_{\D,P}$ exist since $h\in\CaSF'(\dfan^{(s)})$. Now set
$$f_{\D,P}^\tot=f_{\D,P}\cdot \left(\frac{ty_0-y_1}{sy_0-y_1}\right)^{\nu_s(f_{\D,P})}\in \CC(\Y),$$ where $\nu_s$ is the valuation in the point $s$.

\begin{prop}\label{prop:totaldivisora}
        With respect to the open covering $\cX=\bigcup U_{\D,P}^\tot$, the functions $$f_{\D,P}^\tot\cdot \chi^{u_{\D,P}}\in \CC(\cX)$$ define an invariant Cartier divisor on $\cX$ which we denote by $D_h^\tot$.
The restriction of $D_h^\tot$ to $\cX_s$ is $D_h$.
\end{prop}

Before proving the proposition, we illustrate it with an example: 
\begin{ex}[Deforming a divisor on ${C(dP_6)}$]
  Consider the p-divisor $\D=\Delta\otimes \{0\}+\emptyset \otimes \infty$ on $\PP^1$, where $\Delta$ is the hexagon from Figure \ref{fig:md}. The Minkowski decomposition  $\Delta=\Delta^0+\Delta^1$ pictured in this figure gives rise to the p-divisor 
  $$\D^\tot=\Delta^0\otimes(\{0\}\times \base)+\Delta^1\otimes V(ty_0-y_1).$$ Consider some $s\in\base^*$.
 Let $h\in\CaSF(\D^{(s)})$ be given by $h_{|\Delta^i}\equiv1$ for $i=0,1$. This is clearly in $\CaSF'(\D^{(s)})$. The support function $h$ corresponds to the principal divisor $D_h$ defined by the function
 $$f=\frac{y_1(sy_0-y_1)}{y_0^2}.$$

In this example, $\mcP=\emptyset$, so the open cover of $X(\D^\tot)$ just consists of $X(\D^\tot)$.
Since $\nu_s(f)=1$, the divisor $D_h$ lifts to the principal divisor on $X(\D^\tot)$ defined by 
 $$f^\tot=\frac{y_1(ty_0-y_1)}{y_0^2}.$$
 The restriction of this divisor to the special fiber is just the principal divisor defined by $(y_1/y_0)^2$.
\end{ex}

\begin{proof}[Proof of Proposition \ref{prop:totaldivisora}]
        Consider $\D,\D'\in \dfan$ and $P,P'\in \mathcal{P}\cup\{0\}$. It is sufficient to show $$\frac{f_{\D,P}^\tot}{f_{\D',P'}^\tot}\cdot \chi^{u_{\D,P}-u_{\D',P'}}\in H^0\Big(U_{\D,P}^\tot\cap  U_{\D',P'}^\tot, \CO_{\cX}\Big).$$
        Setting $\tilde{u}=u_{\D,P}-u_{\D',P'}$, this is equivalent to showing
\begin{equation*}
        g:=\frac{f_{\D,P}}{f_{\D',P'}}\cdot \left(\frac{ty_0-y_1}{sy_0-y_1}\right)^{\nu_s(f_{\D,P}/f_{\D',P'})}\in H^0\Big(\Y_{\D,P}\cap \Y_{\D',P'},\D^\tot\cap\D'^\tot\big(\tilde{u}\big)\Big),
\end{equation*}
where $\Y_{D,P}$ is the locus of the polyhedral divisor used to define $U_{\D,P}^\tot$.
 This in turn is the same as showing that
\begin{equation}\label{eqn:divval}
        \nu_{D}(g) \geq -(\D^\tot\cap\D'^\tot)_{D}\big(\tilde{u}\big)
\end{equation}
for all prime divisors $D$ contained in $\Y_{\D,P}\cap \Y_{\D',P'}$, where $\nu_D$ is the corresponding valuation. One immediately sees that this is automatically fulfilled unless $D$ is of the form  $V(ty_0-y_1)$  or $Q\times\base$ for some $Q\in\PP^1$, since both sides of the above inequality will be $0$.

Now for $Q\in \PP^1\setminus\{s\}$, $\nu_{Q\times\base}(g)=\nu_Q(f_{\D,P}/f_{\D',P'})$. Furthermore, $\nu_{s\times\base}(g)=0$ and $\nu_{V(ty_0-y_1)}(g)=\nu_s(f_{\D,P}/f_{\D',P'})$. On the other hand, we have
\begin{align*}
        (\D^\tot\cap\D'^\tot)_{Q\times\base}\big(\tilde{u}\big)&=(\D^{(s)}\cap\D'^{(s)})_{Q}\big(\tilde{u}\big);\\
                (\D^\tot\cap\D'^\tot)_{s\times\base}\big(\tilde{u}\big)&=0;\\
                (\D^\tot\cap\D'^\tot)_{V(ty_0-y_1)}\big(\tilde{u}\big)&=(\D^{(s)}\cap\D'^{(s)})_{s}\big(\tilde{u}\big).\\
\end{align*}

Now, since the functions $f_{\D,P}\chi^{u_{\D,P}}$ define a Cartier divisor on $\cX_s$, we have
\begin{equation*}
        \frac{f_{\D,P}}{f_{\D',P'}}\in H^0\Big(Y_{\D,P}\cap Y_{\D',P'},\D^{(s)}\cap\D'^{(s)}\big({\tilde{u}}\big)\Big)
\end{equation*}
where $Y_{D,P}$ is defined similarly to $\Y_{D,P}$. Consequently,
\begin{equation}
        \nu_{Q}\left(\frac{f_{\D,P}}{f_{\D',P'}}\right) \geq -(\D^{(s)}\cap\D'^{(s)})_{Q}\big(\tilde{u}\big)
\end{equation}
for $Q\in Y_{D,P}$ and inequality \eqref{eqn:divval} follows for the required divisors.

The fact that $D_h^\tot$ restricts to $D_h$ on $\cX_s$ follows from the easy observation that the functions $f_{\D,P}^\tot$ restrict to the functions $f_{\D,P}$.
\end{proof}

Having checked that $D_h^\tot$ is indeed a Cartier divisor of $\cX$, we now want to describe its restrictions to the fiber  $\cX_0$. This restriction $(D_h^\tot)_{0}$ will be $T$-invariant, and should thus correspond to some support function $h^{(0)}\in\CaSF(\dfan)$. 

\begin{definition}\label{def:hzero}
Given $h\in\CaSF'(\dfan^{(s)})$, define $h^{(0)}\in \CaSF(\dfan)$ as follows:
\begin{itemize}
        \item For $P\in \PP^1\setminus \{0,s\}$ set $h_P^{(0)}=h_P$;
        \item Set $h_s^{(0)}=h^\lin$;
        \item Finally, set
$$
h_0^{(0)}(v)=h_0(v_0)+h_S(v_1),
$$
where if $v\in\Delta$ for some $\Delta\in\dfan_0$, we take any $v_0\in\Delta^0$, $v_1\in\Delta^1$  such that $v_0+ v_1=v$.
\end{itemize}
Note that the requirement $h\in \CaSF'(\dfan^{(s)})$ ensures that $h_0^{(0)}(v)$ does not depend on the choice of such $v_0$ and $v_1$. One easily checks that $h^{(0)}$ is in fact an element of $\CaSF(\dfan)$.
\end{definition}

\begin{prop}\label{prop:restrict}
Consider $h\in\CaSF'(\dfan^{(s)})$.
 The restriction of $D_h^\tot$ to the special fiber $\cX_0$  is equal to
        $$(D_h^\tot)_{0}=D_{h^{(0)}}.$$
\end{prop}
\begin{proof}
        A straightforward calculation shows that the restrictions of the functions $f_{\D,P}^\tot\cdot \chi^{u_{\D,P}}$ to the fiber $\cX_{0}$ are exactly those determined by $h^{(0)}$.
\end{proof}

In light of the two above propositions, we define a map $\pim:\TCDiv'(\cX_s)\to \TCDiv(\cX_0)$ by sending $D_{h}$ to $D_{h^{(0)}}$.  It is clear from construction that $\pim$ is a group homomorphism sending principal divisors to principal divisors, with kernel contained in the set of principal divisors. Thus, $\pim$ always descends to an injective map $\bpim:\Pic'(\cX_s)\hookrightarrow\Pic(\cX_0)$, where $\Pic'$ is the image of $\CaSF'$ modulo linear equivalence.
In the following section, we shall further study the properties of $\pim$ and $\bpim$. Before doing so however, we express the map $\pim$ in terms of invariant prime divisors, and consider an example.
As in Proposition \ref{prop:weil}, for any $v\in N_\QQ$, let $\mu(v)$ be the smallest natural number such that $\mu(v)v\in N$.
\begin{prop}
Consider some Cartier divisor $D\in\TCDiv'(\cX_s)$, which we can write as
$$
D=\sum_{\rho\in\xray(\dfan)} a_\rho D_\rho+\sum_{\substack{P\in\PP^1\\v\in\xvert_P(\dfan^{(s)}) }}b_{P,v} D_{P,v}.
$$
Then 
$$
\pim(D)=\sum_{\rho\in\xray(\dfan)} a_\rho D_\rho+\sum_{\substack{P\in\PP^1\setminus\{0,s\}\\v\in\xvert_P(\dfan) }}b_{P,v} D_{P,v}+
\sum_{v\in\xvert_0(\dfan)} \mu(v)\left(\frac{b_{0,v^0}}{\mu(v^0)}+\frac{b_{s,v^1}}{\mu(v^1)}\right)D_{0,v}
$$
where the $D_\rho$, $D_{P,v}$ now denote invariant prime divisors on $\cX_0$, and $v=v^0+v^1$ is the Minkowski decomposition of any vertex $v$ of $\dfan_0$.
\end{prop}
\begin{proof}
The claim follows from Proposition \ref{prop:restrict} coupled with Proposition \ref{prop:weil} and a straightforward calculation.
\end{proof}

\begin{ex}[A deformation of $\overline{C(dP_6)}$]\label{ex:dpdiv}
	We return to Example \ref{ex:dpdef} and consider the deformation $\pi$ of $\cX_0=\overline{C(dP_6)}=X(\dfan)$ to $\cX_s=\PP^1\times\PP^1\times\PP^1$. We first observe that  $\TCDiv'(\cX_s)\cong \ZZ^3\times \CDiv^0(\mathbb{P}^1)$, where $\CDiv^0(\mathbb{P}^1)$ consists of degree $0$ divisors on $\PP^1$. Indeed, for $a_1,a_2,a_3\in\ZZ$, let $h[a_1,a_2,a_3]\in \CaSF'(\dfan^{(s)})$ be the support function taking respective values $-a_1,-a_2,-a_3$ on the vertices $(0,0),(0,1),(-1,0)$ of $\dfan_0^{(s)}$,  taking respective values $0,a_2-a_1,a_3-a_1$ on the vertices $(0,0),(0,-1),(1,0)$ of $\dfan_s^{(s)}$ and taking value $0$ on all other vertices. Note that this completely determines $h^{(s)}[a_1,a_2,a_3]$. It is then obvious that any element of $\CaSF'(\dfan^{(s)})$ can be written uniquely as $h^{(s)}[a_1,a_2,a_3]+D$ for some $a_1,a_2,a_3\in \ZZ$ and $D\in \CDiv^0(\mathbb{P}^1)$. This gives the above isomorphism.
	
	On the other hand, we also have that $\TCDiv(\cX_0)\cong \ZZ^3\times \CDiv^0(\mathbb{P}^1)$. Indeed, for $a_1,a_2,a_3\in\ZZ$, let $h^{(0)}[a_1,a_2,a_3]\in \CaSF(\dfan)$ be the support function taking respective values $-a_1,-a_2,-a_3$ at $(0,0),(0,1),(-1,0)$ of $\dfan_0^{(0)}$ and with value $0$ on the vertex $0$ of all other slices. As before, this completely determines $h^{(0)}[a_1,a_2,a_3]$ and as above, any element of $\CaSF(\dfan)$ can be written uniquely as $h^{(0)}[a_1,a_2,a_3]+P$. Now, if we take $h=h^{(s)}[a_1,a_2,a_3]$, then one easily checks that $h^{(0)}=h^{(0)}[a_1,a_2,a_3]$. Factoring out by linear equivalence, we then in fact have
that $\bpim$ is an isomorphism.

Note that in this case, $\Pic'(\cX_s)$ is a strict subgroup of $\Pic(\cX_s)$. Indeed, the former is just $\ZZ$, whereas the latter is $\ZZ^3$. In fact, $\Pic'(\cX_s)$ is generated by half of the anticanonical class. 
\end{ex}

\section{Invariant Families of Divisors II}\label{sec:div2}
We continue the situation of the previous section and now proceed to analyze some properties of the maps $\pim$ and $\bpim$. To begin with, we have the following:

\begin{theorem}\label{thm:euler}
	Let the special fiber $\cX_0$ be complete and consider any $D\in\TCDiv'(\cX_s)$. Then we have 
	\begin{align*}
		 h^i\left(\CO(\pim(D))\right)&\geq h^i\left(\CO(D)\right)\qquad \textrm{for all}\ i\geq 0;\\ 
		 \chi\left(\CO(\pim(D))\right)&=\chi\left(\CO(D)\right).
	\end{align*}
\end{theorem}
\begin{proof}
	Consider $h\in \CaSF'(\dfan^{(s)})$ such that $D_{h}=D$. Then $\CO(D_h^\tot)$ is a line bundle on $\cX$ and thus flat over $\base$, since $\pi$ is flat. Now since $\cX_0$ is complete, we have that $\pi$ is proper by Theorem \ref{thm:tdefprop}. Since the restrictions of $\CO(D_h^\tot)$ to $\cX_s$ and $\cX_0$ are respectively $\CO(D)$ and $\CO(\pim(D))$,
the theorem then follows from cohomology and base change, see for example the corollary in \cite[Section II.5]{mumford:70a}.
\end{proof}

Similarly, if $\cX_0$ is complete, $\pim$ preserves intersection numbers:
\begin{theorem}\label{thm:intersectionnumbers}
Let $\cX_0$ be complete of dimension $n$. Consider invariant divisors $D^1,\ldots,D^n$ in $\TCDiv'(\cX_s)$. Then 
$$
\pim(D^1).\cdots.\pim(D^n)=D^1.\cdots.D^n.
$$
\end{theorem}
\begin{proof}
	By Proposition \ref{prop:totaldivisora}, we can lift the divisor $D^i$ to a divisor $(D^i)^\tot$ on $\cX$. Define $\gamma$ to be the one-cycle class on $\cX$ attained by intersecting the divisors $(D^1)^\tot,\ldots,(D^n)^\tot$. Then $\gamma_{0}$, the restriction of $\gamma$ to $\cX_{0}$, is the intersection of all $\pim(D^i)$. Thus, $\deg(\gamma_{0})$ is the intersection number on the left hand side of the above equation. Likewise, $\deg(\gamma_{s})$ gives the intersection number  on the right, where $\gamma_{s}$ is the restriction of $\gamma$ to $\cX_s$. The theorem then follows from a direct application of Proposition 10.2 in \cite{fulton:98a}.
\end{proof}

We also have that $\pim$ preserves canonical divisors and semiample divisors:
\begin{theorem}\label{thm:canonical}
	If  $K\in \TCDiv'(\cX_s)$ is a canonical divisor on $\cX_s$, then $\pim(K)$ is a canonical divisor on $\cX_0$. If the support of every slice $\dfan_P$ is convex and $D\in\TCDiv'(\cX_s)$ is semiample, then $\pim(D)$ is semiample.
\end{theorem}
\begin{proof}
	If $K\in \TCDiv'(\cX_s)$, we can assume (after possible modification with an invariant principal divisor) that it is of the form stated in \cite[Theorem 3.19]{petersen:08a}.  Coupled with  \cite[Proposition 3.16]{petersen:08a}, we have that $K=D_{h}$, with $h\in\CaSF'(\dfan^{(s)})$ defined as follows:
	\begin{enumerate}
		\item For $P \in Y\setminus\{0\}$ and $v$ a vertex in $\dfan_P^{(s)}$, $h_P(v)=-1+1/\mu(v)$;
		\item For any $v$ a vertex  in $\dfan_0^{(s)}$, $h_0(v)=1+1/\mu(v)$;
		\item The function $h^\lin$ has slope $1$ along every ray of the tailfan of $\dfan^{(s)}$.
	\end{enumerate}
	Indeed, this follows immediately by taking $K_Y=-2\cdot\{0\}$ in  \cite[Theorem 3.19]{petersen:08a}.
On the other hand, $h^{(0)}\in\CaSF(\dfan)$ is the support function defined by:
\begin{enumerate}
	\item For $P \in Y\setminus\{0\}$ and $v$ a vertex  in $\dfan_P$, $h_P^{(0)}(v)=-1+1/\mu(v)$;\label{item:canone}
	\item For $v$ a vertex in $\dfan_0$, $h_0^{(0)}(v)=1+1/\mu(v)$;\label{item:cantwo}
	\item The function $(h^{(0)})^\lin$ has slope $1$ along every ray of the tailfan of $\dfan$.\label{item:canthree}
	\end{enumerate}
	Indeed, \eqref{item:canone} and \eqref{item:canthree} are immediate, and \eqref{item:cantwo} follows from the fact that any vertex $v\in\dfan_0$ is the sum of vertices of $\dfan_0^{(s)}$ and $\dfan_s^{(s)}$, one of which must be a lattice point.	Taking again $K_Y=-2\cdot\{0\}$, we see that $D_{h^{(0)}}$ is also canonical.

Now, any divisor $D_h$ on a complexity-one $T$-variety coming from a divisorial fan whose slices have convex support is semiample if and only if $h$ is concave and satisfies some positivity assumptions, see \cite[Theorem 3.27]{petersen:08a} for the complete case; the general case is similar and follows from \cite[Theorem 3.2]{ilten:11a}.
Suppose now that $D\in\TCDiv'(\dfan^{(s)}$ is semiample, and let $h\in\CaSF(\dfan^{(s)}$ be such that $D=D_h$. Then the concavity and positivity conditions for $h^{(0)}$ follow immediately from those for $h$.
\end{proof}

Finally, we show that if the special fiber $\cX_0$ is semiprojective, i.e. projective over something affine, then $\bpim$ is an isomorphism:
\begin{theorem}
Assume that $\cX_0$ is semiprojective.  Then $\bpim:\Pic'(\cX_s)\to\Pic(\cX_0)$ is surjective, and thus an isomorphism.
\end{theorem}
\begin{proof}
Since $\cX_0$ is semiprojective, any divisor $D$ may be written as the difference of two semiample divisors. Thus, we must only check that for any $f\in\CaSF(\dfan)$ with $D_f$ semiample, there is some $h\in\CaSF'(\dfan^{(s)})$ with $h^{(0)}$ differing from $f$ by a principal divisorial support function.

After correcting by a principal divisorial support function, we may assume that $f_s=f^\lin$. Since $D_f$ is semiample, $f_0$ must be concave, cf. \cite[Corollary 3.3.1]{ilten:11a}. Now, let $\Gamma$ be the polyhedral subdivision
$$
\Gamma=\bigcup_{\Delta\in\dfan_0} \{(v,a)\in \Delta \times \QQ\ |\ a\leq f_0(v)\} 
$$
which has convex support, since $f_0$ is concave. Note that there is a natural projection from $\Gamma$ to $\dfan_0$.
Choose any vertex $e$ of $\dfan_0$ with decomposition $e=e^0+e^1$, where we assume without loss of generality that $e^0$ is a lattice point. Requiring that $(e,f_0(e))^0=(e^0,0)$, the Minkowski decomposition of $\dfan$ determines a Minkowski decomposition of $\Gamma$ compatible with the above projection. By Proposition \ref{prop:mdcomp}, $\Gamma^0$ and $\Gamma^1$ are convex polyhedral subdivisions. Take $h_0$ and $h_s$ to be the concave functions on $\dfan_0^{(s)}$ and $\dfan_s^{(s)}$ determined by the upper faces of $\Gamma^0$ and $\Gamma^1$.
For $P\neq0,s$, set $h_P=f_P$. Then $h$ is a divisorial support function, and if it is integral, it is in $\CaSF'(\dfan^{(s)})$. Furthermore, $h^{(0)}=f$, where we expand Definition \ref{def:hzero} to include non-integral support functions.

Thus, it remains to show that $h_0$ and $h_s$ are integral.
Now for any $\Delta\in \dfan_0$, there is $u\in M$ and $a\in \ZZ$ such that $f_0$ restricted to $\Delta$ is given by $a+\langle\cdot,u\rangle$. By construction there are $a_0,a_s\in \QQ$ such that $a=a_0+a_s$ and $h_0$ (respectively $h_s$) restricted to $\Delta^0$ (or $\Delta^1$) is given by  $a_0+\langle\cdot,u\rangle$ (respectively $a_s+\langle\cdot,u\rangle$). Note that if $a_0\in \ZZ$, then $a_s\in \ZZ$ and vice versa. Thus, for each $\Delta$, we must show that either $a_0\in\ZZ$ or $a_s\in \ZZ$. Note that if $\Delta^0$ or $\Delta^1$ contains a lattice point on which $h_0$ or $h_b$ takes an integral value, then $a_0$ respectively $a_s$ is integral as well. Thus, $h_0$ and $h_b$ are integral on any $\nabla^i$, where $e_0\in \nabla^0$. Now, for any general $\Delta\in\dfan_0$ intersecting such a $\nabla$, $\Delta^i\cap\nabla^i$ must contain a lattice point for either $i=0$ or $i=1$ by the admissibility of the decomposition of $\dfan$.  Thus, $h_0$ and $h_s$ are integral on such $\Delta^i$ as well.  Proceeding by induction using the connectedness of $\dfan_0$ completes the claim.
\end{proof}

\section{Further Remarks and Examples}\label{sec:rems}
We conclude the paper with a number of remarks and examples.

\begin{rem}
The map $\pim$ is never injective. Indeed, the principal divisor $D_h$ on $\cX_s$ with $h_0=1$, $h_s=-1$, and $h_P=0$ for $p\neq 0,s$ maps to the trivial divisor.  Likewise, the map $\pim$ is never surjective.  Indeed, any divisor $D_h$ on $\cX_0$ with $h_s\neq 0$ cannot lie in the image of $\pim$.
\end{rem}

\begin{rem}
It is not at all surprising that the map $\bpim$ is surjective. Indeed, for any rational complexity-one $T$-variety $X$, the cohomology $H^2(X,\CO_X)$ vanishes, cf. \cite[Proposition 38]{tsurvey}. This implies that for every $\mcL\in\Pic(X)$, any first-order deformation $\cX$ of $X$ can be lifted to a first-order deformation of the pair $(X,\mcL)$, see \cite[Theorem 3.3.11]{sernesi:06a}.
\end{rem}

\begin{figure}[h]
  \centering
  \subfigure[Fan for $\F_2$]{\hspace*{2ex}\fanpaab\hspace*{1ex}}
  \subfigure[A Minkowski decomposition]{\paabmink}
  \subfigure[Fan for $\PP^1\times\PP^1$]{\hspace*{3ex}\fanpapa\hspace*{3ex}}
  \caption{Deforming $\F_2$ to $\PP^1\times\PP^1$}
\label{fig:f2tof0}
\end{figure}

Although we have seen that $\pim$ preserves the property of being semiample, it does not in general preserve the property of being ample.
\begin{ex}[A deformation of $\F_2$]
Consider the second Hirzebruch surface $X=\F_2$. This is a toric variety corresponding to the fan $\Sigma$ in $\ZZ^2\otimes \QQ$ pictured in Figure \ref{fig:f2tof0}(a). We can view this as a $T$-variety for $T=\ker r$, where $r$ is the character corresponding to $[0,1]\in(\ZZ^2)^*$. Following the downgrading procedure of \cite[Section 5]{altmann:08a} gives us a divisorial fan $\dfan$ with $X=X(\dfan)$, whose only non-trivial slice is $\dfan_0=\Sigma\cap[r=1]$. This slice, shrunk by a factor of two, is pictured in the middle of Figure \ref{fig:f2tof0}(b). In the same figure, an admissible one-term Minkowski decomposition of $\dfan_0$ is pictured, with the gray line segments denoting which summands belong together. This decomposition gives us a deformation $\pi$, whose general fiber is $\PP^1\times\PP^1$.  Indeed, by reversing the downgrading procedure for the corresponding divisorial fan $\dfan^{(s)}$ $(s\neq 0)$ we get the fan pictured in figure \ref{fig:f2tof0}(c).

Since $\pi$ comes from an essentially locally trivial decomposition, $\bpim$ gives us an isomorphism $\bpim:\Pic(\PP^1\times\PP^1)\to\Pic(\F_2)$.
  However, the image of an ample divisor isn't necessarily ample.  Indeed, the divisor $D_{0,-1}+D_{0,0}$ on $\cX_s\cong\PP^1\times\PP^1$ corresponds to the line bundle $\CO(1,1)$ and is thus ample. This maps via $\pim$ to $D_{0,-1}+D_{0,0}+D_{0,1}$ which is semiample, but not ample.  Indeed, this divisor corresponds to the pullback of $\CO(2)$ under the minimal resolution $\F_2\to\PP(1,1,2)$. 
\end{ex} 

\begin{figure}[h]
  \centering
  \subfigure[Fan for $\PP(1,1,3)$]{\hspace*{2ex}\fanpaac\hspace*{1ex}}
  \subfigure[A Minkowski decomposition]{\paacmink}
  \subfigure[Fan for $\F_1$]{\hspace*{3ex}\fanfa\hspace*{3ex}}
  \caption{Deforming $\PP(1,1,3)$ to $\F_1$}
\label{fig:p113tof1}
\end{figure}

\begin{rem}
	The above example is a special case of the classical deformation of a Hirzebruch surface $\F_m$ to another Hirzebruch surface $\F_n$ with $0\leq n <m$ and $2|(m-n)$, see \cite[Example 1.2.2]{sernesi:06a}. Such deformations can all in fact be constructed as T-deformations, cf. \cite[Section 2.2]{ilten:09a}, and a discussion similar to that above can be used to find ample divisors degenerating to non-ample divisors in this more general situation.
\end{rem}

If the special fiber $\cX_0$ is $\QQ$-Gorenstein and the slice $\dfan_0$ has only integral vertices, then it follows from the description of the canonical divisor used in the proof of Theorem \ref{thm:canonical} that some multiple of the canonical class is contained in $\Pic'(\cX_s)$. However, if $\dfan_0$ has non-integral vertices, then this need not be the case, as we can see in the following example.

\begin{ex}[A smoothing of $\PP(1,1,3)$]
Consider the weighted projective space $X=\PP(1,1,3)$. This is a toric variety corresponding to the fan $\Sigma$ in $\ZZ^2\otimes \QQ$ pictured in Figure \ref{fig:p113tof1}(a). We can view this as a $T$-variety for $T=\ker r$, where $r$ is the character corresponding to $[0,1]\in(\ZZ^2)^*$. Following the dowgrading procedure of \cite[Section 5]{altmann:08a} gives us a divisorial fan $\dfan$ with $X=X(\dfan)$, whose only non-trivial slice is $\dfan_0=\Sigma\cap[r=1]$. This slice, shrunk by a factor of two, is pictured in the middle of Figure \ref{fig:p113tof1}(b). In the same figure, an admissible one-term Minkowski decomposition of $\dfan_0$ is pictured, with the gray line segments denoting which summands belong together. This decomposition gives us a deformation $\pi$, whose general fiber is the Hirzebruch surface $\F_1$.  Indeed, by reversing the downgrading procedure for the corresponding divisorial fan $\dfan^{(s)}$ $(s\neq 0)$ we get the fan pictured in figure \ref{fig:p113tof1}(c).

In this example, the rank of $\Pic(\cX_s)$ is two, whereas that of $\Pic'(\cX_s)$ and $\Pic(\cX_0)$ is one. Now, $\cX_0$ is $\QQ$-Gorenstein with Gorenstein index three; the divisor $5D_{0,-1}+10D_{0,1/2}$ is Cartier and has class thrice that of the anticanonical. However, the inverse of this class under $\bpim$ is given by the class of $5(D_{0,-1}+D_{0,0})$ on $\cX_s$, which is not pluricanonical.  Indeed, the anticanonical class of $\cX_s$ can be represented by $D_{0,-1}+D_{0,0}+D_{s,0}+D_{s,1/2}$, which doesn't lie in $\TCDiv'(\cX_s)$.
\end{ex}

Finally, we say something about deformations of sections. Since (for $\cX_0$ complete) $$h^0(\cX_0,\CO(\pim(D)))\geq h^0(\cX_s,\CO(D),$$ we might expect there to be a natural inclusion $$H^0(\cX_s,\CO(D))\hookrightarrow H^0(\cX_0,\CO(\pim(D))).$$ It is indeed possible to construct such an inclusion, but it is not canonical. For any function $f\in\CC(\PP^1)$, set
$$
f^{(0)}=f\cdot\left(\frac{y_1-sy_0}{y_1}\right)^{-\nu_s(f)}.
$$
Note that for any $u\in M$, $f\in H^0(\cX_s,\CO(D))_u$ implies that $f^{(0)}\in H^0(\cX_0,\CO(\pim(D)))_u$.
Now for every $u\in M$, let $\{b_u^i\}$ be a basis of $H^0(\cX_s,\CO(D))_u$ such that $\{(b_u^i)^{(0)}\}$ gives a basis for $H^0(\cX_0,\CO(\pim(D)))_u$; we leave it to the reader to check that such $\{(b_u^i)^{(0)}\}$ exist.
This then induces the desired inclusion (by sending $b_u^i\cdot\chi^u$ to $(b_u^i)^{(0)}\cdot \chi^u$) but is dependent on the above choice of basis.

\subsubsection*{Acknowledgements} We would like to thank Jos\'e Ignacio Burgos Gil and Martin Sombra for helpful comments.

\bibliographystyle{alpha}
\bibliography{famdiv}

\end{document}

%% file: figures.tex
\newcommand{\dpnull}{%
 \psset{unit=0.5cm}
 \begin{pspicture}(-4,-3.2)(4,3.2)%
\pspolygon*[linecolor=gray4](1,1)(1,0)(0,-1)(-1,-1)(-1,0)(0,1)(1,1)
\pspolygon*[linecolor=gray1](3,0)(1,0)(1,1)(3,3)
\pspolygon*[linecolor=gray6](3,3)(1,1)(0,1)(0,3)
\pspolygon*[linecolor=gray2](0,3)(0,1)(-1,0)(-3,0)(-3,3)
\pspolygon*[linecolor=gray7](-3,0)(-1,0)(-1,-1)(-3,-3)
\pspolygon*[linecolor=gray3](-3,-3)(-1,-1)(0,-1)(0,-3)
\pspolygon*[linecolor=gray5](0,-3)(0,-1)(1,0)(3,0)(3,-3)
   \psgrid[gridwidth=0.3pt,griddots=5,subgriddiv=1,gridlabels=5pt](-3,-3)(3,3)
 \psset{linewidth=1pt}%
\psline(3,3)(1,1)
\psline(3,0)(1,0)
\psline(0,-3)(0,-1)
\psline(-3,-3)(-1,-1)
\psline(-3,0)(-1,0)
\psline(0,3)(0,1)
\psline(1,1)(1,0)(0,-1)(-1,-1)(-1,0)(0,1)(1,1)
\end{pspicture}}

\newcommand{\dpnullia}{%
 \psset{unit=0.5cm}
 \begin{pspicture}(-4,-3.2)(4,3.2)%
\pspolygon*[linecolor=gray4](-1,0)(0,1)(0,0)
\pspolygon*[linecolor=gray1](3,0)(0,0)(0,1)(2,3)(3,3)
\pspolygon*[linecolor=gray6](2,3)(0,1)(0,3)
\pspolygon*[linecolor=gray2](0,3)(0,1)(-1,0)(-3,0)(-3,3)
\pspolygon*[linecolor=gray7](-3,0)(-1,0)(-3,-2)
\pspolygon*[linecolor=gray3](-3,-2)(-1,0)(0,0)(0,-3)(-3,-3)
\pspolygon*[linecolor=gray5](0,-3)(0,0)(3,0)(3,-3)
	 \psgrid[gridwidth=0.3pt,griddots=5,subgriddiv=1,gridlabels=5pt](-3,-3)(3,3)

 \psset{linewidth=1pt}%
\psline(3,0)
\psline(0,-3)
\psline(-3,-2)(2,3)
\psline(-3,0)
\psline(0,3)
\end{pspicture}}

\newcommand{\dpnullib}{%
 \psset{unit=0.5cm}
 \begin{pspicture}(-4,-3.2)(4,3.2)%
\pspolygon*[linecolor=gray4](1,0)(0,-1)(0,0)
\pspolygon*[linecolor=gray1](3,0)(1,0)(3,2)
\pspolygon*[linecolor=gray6](3,2)(1,0)(0,0)(0,3)(3,3)
\pspolygon*[linecolor=gray2](0,3)(0,0)(-3,0)(-3,3)
\pspolygon*[linecolor=gray7](-3,0)(0,0)(0,-1)(-2,-3)(-3,-3)
\pspolygon*[linecolor=gray3](-2,-3)(0,-1)(0,-3)
\pspolygon*[linecolor=gray5](0,-3)(0,-1)(1,0)(3,0)(3,-3)
	 \psgrid[gridwidth=0.3pt,griddots=5,subgriddiv=1,gridlabels=5pt](-3,-3)(3,3)

 \psset{linewidth=1pt}%
\psline(3,0)
\psline(0,-3)
\psline(-3,0)
\psline(0,3)
\psline(-2,-3)(3,2)
\end{pspicture}}

\newcommand{\fanpaab}{
 \psset{unit=1cm}
 \begin{pspicture}(-1.4,-1.4)(1.2,1.5)%
 \psgrid[gridwidth=0.3pt,griddots=5,subgriddiv=1,gridlabels=5pt,gridlabelcolor=white](-1,-1)(1,1)
\psline{<->}(0,-1)(0,1)
\psline{->}(1,1)
\psline{->}(-1,1)
\psline[linecolor=lightgray,linewidth=1.5pt](-1,1)(1,1)
\rput(-1.8,1){\gray{$[r=1]$}}
  \end{pspicture}
}
\newcommand{\fanpapa}{
 \psset{unit=1cm}
 \begin{pspicture}(-1.2,-1.4)(1.05,1.5)%
 \psgrid[gridwidth=0.3pt,griddots=5,subgriddiv=1,gridlabels=5pt,gridlabelcolor=white](-1,-1)(1,1)
\psline{<->}(0,-1)(0,1)
\psline{->}(1,-1)
\psline{->}(-1,1)
  \end{pspicture}
}

\newcommand{\paabmink}{
 \psset{xunit=1cm}
 \begin{pspicture}(-3,-1.4)(1.9,1.5)%
\psline{<->}(-1.3,1)(1.3,1)
\psline{<->}(-1.3,0)(1.3,0)
\psline{<->}(-1.3,-1)(1.3,-1)
\psline[linecolor=lightgray](-1,1)(0,-1)(0,1)(1,-1)
\rput(-2,1){$\{\Delta^0\}$}
\rput(-2,-1){$\{\Delta^1\}$}
\rput(-2.3,0){$\Sigma\cap[r=1]$}
\rput(-1,1.3){$-1$}
\rput(0,1.3){$0$}
\rput(0,-1.3){$0$}
\rput(1,-1.3){$1$}
\end{pspicture}}

\newcommand{\fanpaac}{
 \psset{unit=1cm}
 \begin{pspicture}(-1.4,-2.4)(1.2,2.5)%
 \psgrid[gridwidth=0.3pt,griddots=5,subgriddiv=1,gridlabels=5pt,gridlabelcolor=white](-1,-1)(1,2)
\psline{->}(0,-1)
\psline{->}(1,2)
\psline{->}(-1,1)
\psline[linecolor=lightgray,linewidth=1.5pt](-1,1)(1,1)
\rput(-1.8,1){\gray{$[r=1]$}}
  \end{pspicture}
}
\newcommand{\fanfa}{
 \psset{unit=1cm}
 \begin{pspicture}(-1.2,-2.4)(1.05,1.5)%
 \psgrid[gridwidth=0.3pt,griddots=5,subgriddiv=1,gridlabels=5pt,gridlabelcolor=white](-1,-2)(1,1)
\psline{<->}(0,-1)(0,1)
\psline{->}(1,-2)
\psline{->}(-1,1)
  \end{pspicture}
}

\newcommand{\paacmink}{
 \psset{xunit=1cm}
 \begin{pspicture}(-3,-2.4)(1.9,1.5)%
\psline{<->}(-1.3,1)(1.3,1)
\psline{<->}(-1.3,0)(1.3,0)
\psline{<->}(-1.3,-1)(1.3,-1)
\psline[linecolor=lightgray](-1,1)(0,-1)
\psline[linecolor=lightgray](0,1)(.5,-1)
\rput(-2,1){$\{\Delta^0\}$}
\rput(-2,-1){$\{\Delta^1\}$}
\rput(-2.3,0){$\Sigma\cap[r=1]$}
\rput(-1,1.3){$-1$}
\rput(0,1.3){$0$}
\rput(0,-1.3){$0$}
\rput(.5,-1.3){$\frac{1}{2}$}
\end{pspicture}}